# REMARKS ON THE DERIVATION OF THE VIRIAL IDENTITY FOR NONLINEAR SCHRÖDINGER EQUATIONS

TOMOYUKI IKEDA, SHUJI MACHIHARA, HAYATO MIYAZAKI, AND TOHRU OZAWA

ABSTRACT. We revisit the derivation of the virial identity for nonlinear Schrödinger equations. In [3, 10], several conservation laws, such as for the charge and the energy, were derived without constructing a sequence of approximate solutions. Their approach involves additional properties of solutions due to Strichartz' estimate. In this paper, we derive the virial identity without constructing the sequence of approximate solutions or employing a regularizing argument for weights, by exploiting the properties of solutions.

## 1. INTRODUCTION

In this paper, we consider the Cauthy problem for the nonlinear Schrödinger equation

$$\text{(NLS)} \quad \begin{cases} i\partial_t u + \dfrac{1}{2}\Delta u = f(u), & t \in \mathbb{R},\ x \in \mathbb{R}^n, \\ u(0,x) = \phi(x), & x \in \mathbb{R}^n, \end{cases}$$

where $u = u(t,x)$ is a complex valued function, $\phi \colon \mathbb{R}^n \to \mathbb{C}$ is a given function, and $f$ is a gauge invariant nonlinearity. A typical nonlinearity in our mind is the form of $f(u) = \lambda|u|^{p-1}u$ with $\lambda \in \mathbb{R}$ and $1 < p < \infty$. We assume the following conditions throughout this paper:

(A1) $f \in C^1(\mathbb{C}, \mathbb{C})$, $f(0) = 0$ and for some $1 < p < 1 + 4/(n-2)_+$, $f$ satisfies

$$|f'(z)| := \max\left(\left|\frac{\partial f}{\partial z}\right|, \left|\frac{\partial f}{\partial \bar{z}}\right|\right) \le C(1 + |z|^{p-1})$$

for any $z \in \mathbb{C}$.

(A2) $\operatorname{Im}(\bar{z} f(z)) = 0$ for all $z \in \mathbb{C}$.

(A3) There exists $V \in C^1(\mathbb{C}, \mathbb{R})$ such that $V(0) = 0$ and $f(z) = \partial V/\partial \bar{z}$.

If $f$ satisfies (A2) and (A3), then $V(z) = V(|z|)$ and $f(e^{i\theta}z) = e^{i\theta}f(z)$ hold for any $z \in \mathbb{C}$ and $\theta \in \mathbb{R}$. Hence, we simply write $V'(z) = \frac{dV}{dr}(r)\big|_{r=|z|}$. The aim of this paper is to revisit the derivation of the virial identity for (NLS)

$$\begin{aligned}\|xu(t)\|_{L^2}^2 = \|x\phi\|_{L^2}^2 - 2t\operatorname{Im}\int_{\mathbb{R}^n} x\phi(x)\cdot\overline{\nabla\phi(x)}\,dx + 2t^2 E(\phi)\\ - 2\int_0^t\int_0^s\int_{\mathbb{R}^n} W(u(\tau))\,dxd\tau ds,\end{aligned} \quad (1.1)$$

where $W(u) = (n+2)V(u) - nV'(u)|u|/2$ and $E(\phi)$ is the energy defined in (1.2). Glassey [6] proved the finite time blow-up of solutions to (NLS) by using the identity (1.1). The identity (1.1) is currently known as a fundamental tool to analyze the behavior of the solution.

---

The solution formally possesses the conserved charge $I(u) \coloneqq \|u\|_{L^2}^2$, the conserved energy, and the conserved momentum

$$(1.2) \qquad E(u) \coloneqq \frac{1}{2}\|\nabla u\|_{L^2}^2 + \int_{\mathbb{R}^n} V(u)\,dx, \qquad P(u) \coloneqq \operatorname{Im} \int_{\mathbb{R}^n} u\overline{\nabla u}\,dx.$$

The charge and the energy play an important role to get a priori bound in time in the well-posedness theory of (NLS) in $L^2$ or $H^1$, respectively (cf. [7, 12]). The conserved momentum is applied to blow-up analysis of the solutions (cf. [9]). The solution also satisfies the pseudo-conformal conservation law

$$(1.3) \qquad \|J(t)u(t)\|_{L^2}^2 + 2t^2 \int_{\mathbb{R}^n} V(u(t))\,dx = \|x\phi\|_{L^2}^2 + 2\int_0^t \left(s \int_{\mathbb{R}^n} W(u(s))\,dx\right) ds,$$

where $J(t) = x + it\nabla$. For instance, (1.3) is applied to get a time decay estimate of solutions to (NLS) (cf. [1, Chapter 7]). The last author [10] gave a proof for the conservation law of the charge and the energy by using additional properties of the solution due to Strichartz' estimate, without constructing a sequence of approximate solutions. Fujiwara and the third author [3] proved the conservation of the momentum and the pseudo-conformal conservation law in the same way. To explain the significance of previous works, we describe the standard procedure to derive the energy conservation law. One denotes the scalar product in $L^2$ by

$$(f,g) \coloneqq \int_{\mathbb{R}^n} f(x)\overline{g(x)}\,dx$$

for all $f, g \in L^2(\mathbb{R}^n)$. Formally, the real part of the scalar product between the equation (NLS) and $-\partial_t u$ can be computed as follows:

$$(1.4) \qquad \begin{aligned} 0 &= 2\operatorname{Re}\left(i\partial_t u + \frac{1}{2}\Delta u - f(u), -\partial_t u\right) \\ &= -\operatorname{Re}(\Delta u, \partial_t u) + 2\operatorname{Re}(f(u), \partial_t u) = \frac{d}{dt}E(u(t)). \end{aligned}$$

Other conservation laws can also be formally derived in a similar manner. Considering $H^1$-solutions, which are the natural choice for the solution in view of the energy, each term in (1.4) does not make sense, because $\Delta u, \partial_t u \notin L^2$ for instance. To justify (1.4), at least $H^2$-solutions are required. To this end, we construct a sequence of smooth approximate solutions, by approximating the initial data, or considering the a sequence of regularized equations, which are compatible with the argument for well-posedness theory and formal computations such as (1.4). Note that the both approaches involve a limiting procedure on approximate solutions. The aforementioned studies demonstrate that various conservation laws can be derived, without constructing such a sequence of approximate solutions, by using the properties of the solution to the integral equation corresponding to (NLS). In order to prove the well-posedness of (NLS), an auxiliary space due to Strichartz' estimate is needed as a solution space. By investigating what functions with $u$ belong to this auxiliary space and its dual space and other associated spaces, it is possible to derive conservation laws without constructing the sequence of approximate solutions.

In this paper, by exploiting additional properties of the solution as stated above, we derive the virial identity (1.1). Formally, (1.1) is given as follows: Taking the real part of the scalar product between the equation (NLS) and $i|x|^2 u$, we have

$$(1.5) \qquad 0 = 2\operatorname{Re}\left(i\partial_t u + \frac{1}{2}\Delta u - f(u), i|x|^2 u\right) = \frac{d}{dt}\|xu(t)\|_{L^2}^2 + 2\operatorname{Im}(xu(t), \nabla u(t)).$$

Furthermore, taking the real part of the scalar product between the equation (NLS) and $-(2x \cdot \nabla u + nu)$, one gets

$$
\begin{aligned}
0 &= \mathrm{Re}\left(i\partial_t u + \frac{1}{2}\Delta u - f(u), -(2x \cdot \nabla u + nu)\right) \\
&= \frac{d}{dt}\mathrm{Im}\,(xu(t), \nabla u(t)) + \|\nabla u(t)\|_{L^2}^2 - \int_{\mathbb{R}^n} W(u(t))\,dx.
\end{aligned}
\tag{1.6}
$$

Combining (1.5) with (1.6), together with the conservation of the energy, (1.1) can be obtained. To justify the computation of (1.5) and (1.6), we need to employ not only the construction of the sequence of approximate solutions, but also a regularizing argument with the factor $e^{-\varepsilon|x|^2}$ (see [1, §6.5]). In this paper, it is exhibited that by exploiting additional properties of the solution, (1.1) can be derived without constructing the sequence of approximate solutions or employing the regularizing factor $e^{-\varepsilon|x|^2}$ to ensure spatial integrability. To the best of our knowledge, our method is the first rigorous derivation of the virial identity without using the regularizing factor $e^{-\varepsilon|x|^2}$.

Based on the Duhamel principle, we give the definition of solutions to (NLS) as follows:

**Definition 1** (Solution). *Let $I \ni 0$ be an interval and fix $t_0 \in I$. Take an admissible pair $(q,r) = (4(p+1)/(n(p-1)), p+1)$. We say a function $u\colon I \times \mathbb{R}^n \to \mathbb{C}^n$ is a $H^1$-solution to (NLS) on $I$ if $u \in C(I; H^1(\mathbb{R}^n)) \cap L^q(I; W^{1,r}(\mathbb{R}^n))$ satisfies*

$$
u(t) = U(t-t_0)u(t_0) - i\int_{t_0}^t U(t-s)f(u(s))\,ds
\tag{1.7}
$$

*in $H^1(\mathbb{R}^n)$ for any $t \in I$, where $U(t) = e^{\frac{it\Delta}{2}}$. In addition, $u$ is said to be a $\Sigma$-solution if $u$ satisfies $J(\cdot)u \in C(I; L^2(\mathbb{R}^n)) \cap L^q(I; L^r(\mathbb{R}^n))$.*

Let us first recall Strichartz' estimates. A pair $(\rho, \gamma)$ is said to be admissible if

$$
\rho, \gamma \in [2, \infty], \quad \frac{2}{\rho} = d\left(\frac{1}{2} - \frac{1}{\gamma}\right), \quad (d, \rho, \gamma) \neq (2, 2, \infty).
$$

**Lemma 2** (Strichartz' estimate, e.g., [5, 8, 11, 13]). *Let $(\rho, \gamma)$ and $(\widetilde{\rho}, \widetilde{\gamma})$ be admissible pairs. For any interval $I \ni 0$,*

$$
\|U(\cdot)f\|_{L^\rho(\mathbb{R}; L^\gamma(\mathbb{R}^n))} \leqslant C_0\|f\|_{L^2(\mathbb{R}^n)},
$$

$$
\left\|\int_0^{\cdot} U(\cdot - s)F(s)\,ds\right\|_{L^\rho(I; L^\gamma(\mathbb{R}^n))} \leqslant C_0\|F\|_{L^{\widetilde{\rho}'}(I; L^{\widetilde{\gamma}'}(\mathbb{R}^n))},
$$

*where $C_0 > 0$ is a certain constant not depending on $I$, and $\eta'$ is the dual exponent of $\eta \geqslant 1$.*

*Remark* 3. Thanks to Strichartz' estimate, if $u$ is a $H^1$-solution, then $u \in L^\rho(I; W^{1,\gamma})$ for any admissible pair $(\rho, \gamma)$. As for the well-posedness results for (NLS), see [7] for $H^1$-solutions, [2, 4] for $\Sigma$-solutions, for instance. We also refer to [1] as a lecture note.

1.1. **Main Result.** Recall that we assume (A1), (A2) and (A3) throughout this paper. In [3], using the properties of the solution to (1.7), the pseudo-conformal conservation law was derived.

**Lemma 4** ([3]). *Let $u$ be a $\Sigma$-solution to (NLS) on $[-T, T]$. Then (1.3) holds for any $t \in [-T, T]$.*

The following identity is corresponding to the identity (1.6):

**Proposition 5.** *Let $u$ be a $\Sigma$-solution to* (NLS) *on $[-T, T]$. Then it holds that*

$$\text{Im}\,(xu(t), \nabla u(t)) = \text{Im}\,(x\phi, \nabla \phi) - 2tE(\phi) + \int_0^t \left( \int_{\mathbb{R}^n} W(u(s))\,dx \right) ds \tag{1.8}$$

*for any $t \in [-T, T]$, where $W(u) = (n+2)V(u) - nV'(u)|u|/2$.*

Combining Lemma 4 with Proposition 5, we establish the virial identity without using the regularizing argument involving $e^{-\varepsilon|x|^2}$.

**Corollary 6** (Virial identity)**.** *Let $u$ be a $\Sigma$-solution to* (NLS) *on $[-T, T]$. Then the virial identity* (1.1) *holds for any $t \in [-T, T]$.*

Let us explain the reason why we do not need to use the regularizing factor $e^{-\varepsilon|x|^2}$ to ensure spatial integrability. The pseudo-conformal conservation law (1.3) is formally equivalent to

$$t\int_{\mathbb{R}^n} W(u(t))\,dx = \frac{d}{dt}\left( \|J(t)u(t)\|_{L^2}^2 + 2t^2 \int_{\mathbb{R}^n} V(u(t))\,dx \right).$$

Collecting the conservation of the energy and (1.5), we formally calculate

$$t\int_{\mathbb{R}^n} W(u(t))\,dx = \frac{d}{dt}\left( \|xu(t)\|_{L^2}^2 + 2t\,\text{Im}\,(xu(t), \nabla u(t)) + 2t^2 E(\phi) \right)$$

$$= t\left( 2E(\phi) + \frac{d}{dt}\,\text{Im}\,(xu(t), \nabla u(t)) \right),$$

which implies (1.8). Also, integrating the first equality in the above over $[0, t]$, we have

$$\int_0^t \left( s \int_{\mathbb{R}^n} W(u(s))\,dx \right) ds = \|xu(t)\|_{L^2}^2 + 2t\,\text{Im}\,(xu(t), \nabla u(t)) + 2t^2 E(\phi) - \|x\phi\|_{L^2}.$$

Plugging (1.8) into the above, (1.1) is established. We translate this procedure into the approach that utilizes the properties of the solutions provided by Strichartz' estimate. The operator $J(t) = x + it\nabla = U(t)xU(-t)$, which is connected to the Galilean symmetry of the equation (NLS), has a property

$$J(t) = \mathcal{M}(t)it\nabla\mathcal{M}(-t), \quad \mathcal{M}(t) := e^{\frac{i|x|^2}{2t}}, \quad t \neq 0. \tag{1.9}$$

With the property in place, we have further additional properties of the solutions. The additional properties enable us to derive the virial identity without using the regularizing factor $e^{-\varepsilon|x|^2}$.

**Notations.** For any $p \geqslant 1$, $L^p = L^p(\mathbb{R}^n)$ denotes the usual Lebesgue space on $\mathbb{R}^n$. $\eta' \geqslant 1$ is the dual exponent of $\eta \geqslant 1$ defined by $1/\eta + 1/\eta' = 1$. Set $\langle a \rangle = (1 + |a|^2)^{1/2}$ for $a \in \mathbb{C}$ or $\mathbb{R}^n$. $\mathcal{F}[u] = \widehat{u}$ is the usual Fourier transform of a function $u$ on $\mathbb{R}^n$, and $\mathcal{F}^{-1}[u] = \check{u}$ is its inverse. For $m \in \mathbb{Z}$ and $r \in (1, \infty)$, the standard Sobolev space on $\mathbb{R}^n$ is defined by $W^{m,r} = W^{m,r}(\mathbb{R}^n) := \{u \in \mathscr{S}'(\mathbb{R}^n) \mid \|u\|_{W^{m,r}} := \left\|(1-\Delta)^{\frac{m}{2}}u\right\|_{L^r} < \infty\}$, where $\mathscr{S}'(\mathbb{R}^n)$ stands for the space of tempered distributions on $\mathbb{R}^n$, and $(1-\Delta)^{\frac{m}{2}} = \mathcal{F}^{-1}\langle \xi \rangle^m \mathcal{F}$ is the Bessel potential operator. We denote $H^m = W^{m,2}$ for $m \in \mathbb{Z}$. For $T > 0$ and Banach space $X$, we abbreviate $L^p(-T, T; X) = L_t^p X$. $A \lesssim B$ denotes $A \leqslant CB$ for some constant $C > 0$. The notation $\langle \cdot, \cdot \rangle_{B' \times B}$ represents the duality pairing on $B \times B'$, where $B$ is a Banach space such that $B \hookrightarrow L^2 \hookrightarrow B'$ with dense embeddings. More precisely, if $B = W^{m,r'}(\mathbb{R}^n)$ with $m = 0, 1$ and $1 < r < \infty$, then $B' = W^{-m,r}(\mathbb{R}^n)$ and $\langle \cdot, \cdot \rangle_{B \times B'}$ is defined by

$$\langle f, g \rangle_{W^{-m,r} \times W^{m,r'}} = \int_{\mathbb{R}^n} (1-\Delta)^{-\frac{m}{2}} f(x) \cdot (1-\Delta)^{\frac{m}{2}} g(x)\,dx.$$

In particular, if $r = 2$, then $B = H^m(\mathbb{R}^n)$, $B' = H^{-m}(\mathbb{R}^n)$, and $\langle f, g \rangle_{H^{-m} \times H^m} = (\hat{f}, \overline{\hat{g}})$, where $(\cdot, \cdot)$ is the scalar product in $L^2(\mathbb{R}^n)$. Also, for two Banach spaces $X$ and $Y$, we denote by $X + Y$ the Banach space consisting of all elements $f$ such that $f = f_1 + f_2$ with $f_1 \in X$ and $f_2 \in Y$, equipped with the norm $\|f\|_{X+Y} = \inf_{f=f_1+f_2} (\|f_1\|_X + \|f_2\|_Y)$. Further, the Banach space $X \cap Y$ is equipped with the norm $\|f\|_{X \cap Y} = \|f\|_X + \|f\|_Y$. Since $(X + Y)' = X' \cap Y'$, the duality pairing on $(X' \cap Y') \times (X + Y)$ is defined by

$$\langle f, g \rangle_{(X' \cap Y') \times (X+Y)} = \langle f, g_1 \rangle_{X' \times X} + \langle f, g_2 \rangle_{Y' \times Y}$$

for $f \in X' \cap Y'$ and $g = g_1 + g_2 \in X + Y$. Remark that $\langle f, g \rangle_{(X' \cap Y') \times (X+Y)}$ is independent of the decomposition $g = g_1 + g_2$.

*Remark* 7 (Derivation of the conserved momentum). In [3], a solution $u$ in an appropriate space satisfies the following identity:

$$P(u(t)) - P(\phi) = -2 \int_0^t \mathrm{Re} \left\langle \nabla u(\tau), \overline{f(u(\tau))} \right\rangle_{L^2 \cap L^{p+1} \times (L^2 + L^{\frac{p+1}{p}})} d\tau$$

for all $t \in [-T, T]$. They exploit the standard density argument to show the disappearance of the time integral in the right-hand side, and hence $P(u(t)) = P(\phi)$ holds for any $t$. However, the left-hand side of the above identity is pure imaginary and the right-hand side is real, so that $P(u(t)) = P(\phi)$ is immediate. We thus do not need to use the density argument.

## 2. Preliminaries

Throughout this paper, we assume the specific choice of exponents

$$(q, r) = \left( \frac{4(p+1)}{n(p-1)}, p+1 \right)$$

given in Definition 1.

### 2.1. Properties of solutions.
We summarize the properties of $\Sigma$-solutions and $H^1$-solutions to (NLS) used in the proof of Proposition 5. First, by the definition of $\Sigma$-solutions and (1.9), the solutions satisfy

(2.1)
$$f(u) \in L^1(-T, T; H^1) + L^{q'}(-T, T; W^{1,r'}),$$
$$J(\cdot)f(u) \in L^1(-T, T; L^2) + L^{q'}(-T, T; L^{r'})$$

for any $T > 0$. Since $xu = J(t)u - it\nabla u \in L^2$ for any $t$, $\Sigma$-solutions also satisfy

(2.2)
$$xu \in C([-T, T]; L^2) \cap L^q(-T, T; L^r),$$
$$xf(u) \in L^1(-T, T; L^2) + L^{q'}(-T, T; L^{r'}).$$

Further, the additional properties of $H^1$-solutions are given as follows:

(2.3) $\qquad f(u), \Delta u, \partial_t u \in C([-T, T]; H^{-1}) \cap L^q(-T, T; W^{-1,r}).$

Indeed, we have $H^1 \hookrightarrow L^r$ if $p \in (1, 1 + 4/(d-2)_+)$ and hence $L^{r'} \hookrightarrow H^{-1}$. This yields $f(u) \in C([-T, T]; H^{-1})$. Hence, combining this property with the fact that $U(t)$ is a strongly continuous one-parameter unitary group on $H^{-1}$, it follows from the integral equation (1.7) that $u$ is strongly differentiable with respect to $t$ in $H^{-1}$ and $u$ satisfies the equation (NLS) in $H^{-1}$. Further, when $1 < p < 1 + 4/(d-2)_+$, we have $W^{1,r'} \hookrightarrow L^2 \hookrightarrow W^{-1,r}$. By $H^1 \hookrightarrow L^r$ and Hölder's inequality, $f(u) \in L^q(-T, T; W^{-1,r})$ holds. Hence the equation (NLS) yields $\partial_t u \in C([-T, T]; H^{-1}) \cap L^q(-T, T; W^{-1,r})$. Returning to the equation (NLS), we prove the following:

**Lemma 8.** *Let $u$ be a $H^1$-solution to* (NLS). *Then*
$$\operatorname{Im}\left\langle \nabla u, \overline{\nabla f(u)} \right\rangle_{(L^2\cap L^r)\times(L^2+L^{r'})} = 2\operatorname{Re}\left\langle \partial_t u, \overline{f(u)} \right\rangle_{(H^{-1}\cap W^{-1,r})\times(H^1+W^{1,r'})}.$$

Formally, Lemma 8 is obtained by multiplying both sides of (NLS) by $\overline{f(u)}$, integrating over $\mathbb{R}^n$, and taking the imaginary part. For the reader's convenience, we give a rigorous proof in Appendix A, although the proof has been already given in [10].

2.2. **General facts.** We collect some general facts necessary to prove Proposition 5. First, since $f(v)\bar{v} = V'(v)|v|/2$ and $2\operatorname{Re}(f(v)\overline{\nabla v}) = \nabla V(v)$, the following lemma immediately follows from the standard density argument:

**Lemma 9.** *Let $v \in H^1 \cap W^{1,r}$ with $xv \in L^2 \cap L^r$. Then it holds that*
$$\operatorname{Re}\left\langle xv, \overline{\nabla f(v)} \right\rangle_{(L^2\cap L^r)\times(L^2+L^{r'})} - \operatorname{Re}\left\langle \nabla v, \overline{xf(v)} \right\rangle_{(L^2\cap L^r)\times(L^2+L^{r'})}$$
$$= -\frac{n}{2}\int_{\mathbb{R}^n} V'(v)|v|\,dx + n\int_{\mathbb{R}^n} V(v)\,dx.$$

Further, Fubini's theorem and the unitarity of $U(t)$ yield the following:

**Lemma 10.** *Let $\psi \in L^2 \cap L^r$ and $g_j \in L^1(-T,T;L^2) + L^{q'}(-T,T;L^{r'})$ for $j=1,2$. Then, it holds that for any $t \in [-T,T]$,*
$$\left(\psi, \int_0^t U(-s)g_1(s)\,ds\right) = \int_0^t \left\langle U(s)\psi, \overline{g_1(s)} \right\rangle ds,$$

*and*
$$\left(\int_0^t U(-s)g_1(s)\,ds, \int_0^t U(-\theta)g_2(\theta)\,d\theta\right)$$
$$= \int_0^t \left\langle \int_0^\theta U(\theta-s)g_1(s)\,ds, \overline{g_2(\theta)} \right\rangle d\theta + \int_0^t \left\langle g_1(s), \overline{\int_0^s U(s-\theta)g_2(\theta)\,d\theta} \right\rangle ds,$$

*where all of the duality pairings are understood as those on $(L^2 \cap L^r) \times (L^2 + L^{r'})$.*

Formally, we obtain Lemma 10 by switching the order of integration and the unitarity of $U(t)$. A rigorous proof is provided in Appendix A. The last lemma is concerned with the absolute continuity of $\int_{\mathbb{R}^n} V(v(t,x))\,dx$ for a suitable function $v$.

**Lemma 11.** *Let $v \in C([-T,T];H^1) \cap C^1([-T,T];H^{-1})$ with $\nabla v \in L^q(-T,T;L^r)$ and $\partial_t u \in L^q(-T,T;W^{-1,r})$. Then $\int_{\mathbb{R}^n} V(v(\cdot,x))\,dx \in W^{1,1}(-T,T)$. Moreover,*
$$\frac{d}{dt}\left(\int_{\mathbb{R}^n} V(v(t,x))\,dx\right) = 2\operatorname{Re}\left\langle \partial_t v(t), \overline{f(v(t))} \right\rangle_{(H^{-1}\cap W^{-1,r})\times(H^1+W^{1,r'})}$$

*holds in $L^1(-T,T)$.*

Since the proof of Lemma 11 is a bit lengthy, we only give a remark here and postpone the proof to the appendix A. Lemma 11 implies that

(2.4) $$\int_0^t \frac{d}{ds}\int_{\mathbb{R}^n} V(v(s))\,dxds = \int_{\mathbb{R}^n} V(v(t))\,dx - \int_{\mathbb{R}^n} V(v(0))\,dx$$

holds for any $t \in [-T,T]$. Although the fact that $\int_{\mathbb{R}^n} V(v(t,x))\,dx \in W^{1,1}(-T,T)$ was not explicitly mentioned in [3, 10], it provides a rigorous justification for the derivation of the conservation of the energy and the pseudo-conformal conservation law in [3, 10].

## 3. Derivation of the virial identity

### 3.1. Proof of Proposition 5 and Corollary 6.
We shall derive the virial identity in Corollary 6 without the construction of any approximate solutions and the regularizing argument with $e^{-\varepsilon|x|^2}$. The key of the derivation is to involve the pseudo-conformal conservation law (1.3) and the identity (1.8). By means of the conservation of the energy, a direct calculation shows

$$(3.1) \quad \|J(t)u\|_{L^2}^2 + 2t^2 \int_{\mathbb{R}^n} V(u(t))\,dx = \|xu\|_{L^2}^2 + 2t\,\mathrm{Im}\,(xu, \nabla u) + 2t^2 E(\phi).$$

Plugging (1.8) into (3.1), we reach the desired identity from (1.3).

*Proof of Proposition 5.* Applying the operator $J(t) = x + it\nabla$ to the equation (1.7) with $t_0 = 0$, we see from (1.9) that

$$(3.2) \quad J(t)u(t) = U(t)x\phi - i\int_0^t U(t-s)J(s)f(u(s))\,ds.$$

Similarly, it holds that

$$(3.3) \quad \nabla u(t) = U(t)\nabla\phi - i\int_0^t U(t-s)\nabla f(u(s))\,ds.$$

Note that $\mathrm{Im}\,(xu, \nabla u) = \mathrm{Im}\,(J(t)u, \nabla u) - t\|\nabla u\|_{L^2}^2$. Combining (3.2) with (3.3), we see from Lemma 10 that

$$\begin{aligned}
\mathrm{Im}\,(J(t)u, \nabla u) &= \mathrm{Im}\,(x\phi, \nabla\phi) \\
&\quad + \mathrm{Im}\left(x\phi, -i\int_0^t U(-s)\nabla f(u(s))\,ds\right) \\
&\quad + \mathrm{Im}\left(-i\int_0^t U(-s)J(s)f(u(s))\,ds, \nabla\phi\right) \\
&\quad + \mathrm{Im}\left(-i\int_0^t U(-s)J(s)f(u(s))\,ds, -i\int_0^t U(-\theta)\nabla f(u(\theta))\,d\theta\right) \\
&= \mathrm{Im}\,(x\phi, \nabla\phi) \\
&\quad + \mathrm{Re}\int_0^t \left\langle U(\theta)x\phi, \overline{\nabla f(u(\theta))}\right\rangle d\theta \\
&\quad - \mathrm{Re}\int_0^t \left\langle U(s)\nabla\phi, \overline{J(s)f(u(s))}\right\rangle ds \\
&\quad + \mathrm{Re}\int_0^t \left\langle -i\int_0^\theta U(\theta-s)J(s)f(u(s))\,ds, \overline{\nabla f(u(\theta))}\right\rangle d\theta \\
&\quad - \mathrm{Re}\int_0^t \left\langle -i\int_0^s U(s-\theta)\nabla f(u(\theta))\,d\theta, \overline{J(s)f(u(s))}\right\rangle ds,
\end{aligned}$$

which implies

$$(3.4) \quad \begin{aligned}\mathrm{Im}\,(J(t)u, \nabla u) &= \mathrm{Im}\,(x\phi, \nabla\phi) \\ &\quad + \mathrm{Re}\int_0^t \left\langle J(s)u(s), \overline{\nabla f(u(s))}\right\rangle ds - \mathrm{Re}\int_0^t \left\langle \nabla u(s), \overline{J(s)f(u(s))}\right\rangle ds.\end{aligned}$$

Note that all of the time integrals of the duality pairings in the last two equalities make sense as the duality pairing on $(L_t^\infty L^2 \cap L_t^q L^r) \times (L_t^1 L^2 + L_t^{q'} L^{r'})$ with $(q, r) =$

$(4(p+1)/(n(p-1)), p+1)$. It then holds that

$$\operatorname{Re}\left\langle J(s)u, \overline{\nabla f(u)}\right\rangle = \operatorname{Re}\left\langle xu, \overline{\nabla f(u)}\right\rangle - s\operatorname{Im}\left\langle \nabla u, \overline{\nabla f(u)}\right\rangle,$$

$$\operatorname{Re}\left\langle \nabla u, \overline{J(s)f(u)}\right\rangle = \operatorname{Re}\left\langle \nabla u, \overline{xf(u)}\right\rangle + s\operatorname{Im}\left\langle \nabla u, \overline{\nabla f(u)}\right\rangle,$$

where all of terms in the above also make sense as the duality pairing on $(L^2 \cap L^r) \times (L^2 + L^{r'})$, due to (2.1), (2.2) and (2.3). Recalling $J(t) = x + it\nabla$, we see from Lemma 8 and Lemma 9 that

$$\operatorname{Re}\left\langle J(s)u, \overline{\nabla f(u)}\right\rangle - \operatorname{Re}\left\langle \nabla u, \overline{J(s)f(u)}\right\rangle$$
$$= \operatorname{Re}\left\langle xu, \overline{\nabla f(u)}\right\rangle - \operatorname{Re}\left\langle \nabla u, \overline{xf(u)}\right\rangle - 2s\operatorname{Im}\left\langle \nabla u, \overline{\nabla f(u)}\right\rangle$$
$$= -\frac{n}{2}\int_{\mathbb{R}^n} V'(u)|u|\,dx + n\int_{\mathbb{R}^n} V(u)\,dx - 4s\operatorname{Re}\left\langle \partial_t u, \overline{f(u)}\right\rangle_{(H^{-1}\cap W^{-1,r})\times(H^1+W^{1,r'})}.$$

Integrating the above over $[0,t]$, Lemma 11, and the integration by parts imply

$$\operatorname{Re}\int_0^t \left\langle J(s)u(s), \overline{\nabla f(u(s))}\right\rangle ds - \operatorname{Re}\int_0^t \left\langle \nabla u(s), \overline{J(s)f(u(s))}\right\rangle ds$$
$$= -\int_0^t \left(\frac{n}{2}\int_{\mathbb{R}^n} V'(u(s))|u(s)|\,dx - n\int_{\mathbb{R}^n} V(u(s))\,dx\right) ds$$
$$\quad - 2\int_0^t \left(s\frac{d}{ds}\int_{\mathbb{R}^n} V(u(s))\,dx\right) ds$$
$$= -\frac{n}{2}\int_0^t \left(\int_{\mathbb{R}^n} V'(u(s))|u(s)|\,dx\right) ds + (n+2)\int_0^t \left(\int_{\mathbb{R}^n} V(u(s))\,dx\right) ds$$
$$\quad - 2t\int_{\mathbb{R}^n} V(u(t))\,dx.$$

Note that Lemma 11 also implies that $\int_{\mathbb{R}^n} V(u(t))\,dx$ is absolutely continuous in $t$, which allows us to apply the integration by parts. Hence, plugging the above into (3.4) and using the conservation of the energy, we conclude that

$$\operatorname{Im}(xu, \nabla u) = \operatorname{Im}(J(t)u, \nabla u) - t\|\nabla u\|_{L^2}^2$$
$$= \operatorname{Im}(x\phi, \nabla\phi) - 2tE(\phi) + \int_0^t \left(\int_{\mathbb{R}^n} W(u(s))\,dx\right) ds.$$

This is nothing but (1.8), since $W(u) = (n+2)V(u) - nV'(u)|u|/2$. □

*Proof of Corollary 6.* Substituting (1.8) into (3.1), we have

$$\|J(t)u\|_{L^2}^2 + 2t^2\int_{\mathbb{R}^n} V(u(t))\,dx$$
$$= \|xu\|_{L^2}^2 + 2t\operatorname{Im}(x\phi, \nabla\phi) - 2t^2 E(\phi) + 2t\int_0^t \left(\int_{\mathbb{R}^n} W(u(s))\,dx\right) ds.$$

Plugging the above into (1.3), we conclude the desired identity, since $\int_{\mathbb{R}^n} W(u(t))\,dx$ is continuous in $t$, and the integration by parts

$$\int_0^t sg(s)\,ds = \int_0^t s\frac{d}{ds}\int_0^s g(\tau)\,d\tau\,ds = t\int_0^t g(\tau)\,d\tau - \int_0^t \int_0^s g(\tau)\,d\tau\,ds.$$

□

## Appendix A. Proof of Lemmas in section 2

This appendix is devoted to giving the proofs of Lemmas 8, 10 and 11.

*Proof of Lemma 8.* Firstly, we employ Yosida's approximation $(1-\varepsilon\Delta)^{-1}$ as follows:

$$\left\langle \nabla u, \overline{\nabla f(u)} \right\rangle_{(L^2 \cap L^r) \times (L^2 + L^{r'})} = \lim_{\varepsilon \to +0} \mathrm{Im} \left\langle \nabla u, \overline{(1-\varepsilon\Delta)^{-1} \nabla f(u)} \right\rangle_{(L^2 \cap L^r) \times (L^2 + L^{r'})}.$$

Note that the fact $L^{r'} \hookrightarrow H^{-1}$ yields $(1-\varepsilon\Delta)^{-1} \nabla f(u) \in H^1$. Since the equation (NLS) is valid in $H^{-1}$, it follows from Plancherel's theorem that

$$\mathrm{Im} \left\langle \nabla u, \overline{(1-\varepsilon\Delta)^{-1} \nabla f(u)} \right\rangle_{(L^2 \cap L^r) \times (L^2 + L^{r'})} = -\mathrm{Im} \left\langle \Delta u, \overline{(1-\varepsilon\Delta)^{-1} f(u)} \right\rangle_{H^{-1} \times H^1}$$
$$= 2\mathrm{Re} \left\langle \partial_t u, \overline{(1-\varepsilon\Delta)^{-1} f(u)} \right\rangle_{H^{-1} \times H^1},$$

where we use the following fact in the last line:

$$\left\langle f(u), \overline{(1-\varepsilon\Delta)^{-1} f(u)} \right\rangle_{H^{-1} \times H^1} = \int_{\mathbb{R}^n} (1+\varepsilon|\xi|^2)^{-1} |\widehat{f(u)}(\xi)|^2 \, d\xi \in \mathbb{R}.$$

Thanks to $f(u) \in H^1 + W^{1,r'}$ and $\partial_t u \in H^{-1} \cap W^{-1,r}$,

$$\lim_{\varepsilon \to +0} \left\langle \partial_t u, \overline{(1-\varepsilon\Delta)^{-1} f(u)} \right\rangle_{H^{-1} \times H^1} = \left\langle \partial_t u, \overline{f(u)} \right\rangle_{(H^{-1} \cap W^{-1,r}) \times (H^1 + W^{1,r'})}$$

holds. This completes the proof. □

*Proof of Lemma 10.* We only give the proof of the second assertion, since the first one is similar and easier. To justify the calculations in Lemma 10, Yosida's approximation $(1-\varepsilon\Delta)^{-1}$ is involved. By $L^{r'} \hookrightarrow H^{-1}$, one sees from $g_1, g_2 \in L_t^\infty L^2 + L_t^{q'} L^{r'}$ that

$$(1-\varepsilon\Delta)^{-1} g_1, \ (1-\varepsilon\Delta)^{-1} g_2 \in L_t^\infty L^2.$$

This implies that $U(-s)(1-\varepsilon\Delta)^{-1} g_1(s) \overline{U(-\theta)(1-\varepsilon\Delta)^{-1} g_2(\theta)}$ is integrable in $[0,t]^2 \times \mathbb{R}^n$. Combining the unitarity of $U(t)$ in $L^2$ with Fubini's theorem, and then taking the limit as $\varepsilon \to +0$, we prove the second assertion. This completes the proof. □

*Proof of Lemma 11.* First, $\int_{\mathbb{R}^n} V(u(\cdot)) \, dx \in L^1(-T,T)$ immediately follows from the continuity of $\int_{\mathbb{R}^n} V(u(t)) \, dx$ in $t$. We here denote $\widetilde{v}$ by a suitable extension of $v$ in $\mathbb{R} \times \mathbb{R}^n$. Let us prove

$$\lim_{h \to 0} \int_{\mathbb{R}^n} \frac{V(\widetilde{v}(t+h)) - V(\widetilde{v}(t))}{h} \, dx = 2\mathrm{Re} \left\langle \partial_t \widetilde{v}(t), \overline{f(\widetilde{v}(t))} \right\rangle_{(H^1 + W^{1,r'}) \times (H^{-1} \cap W^{-1,r})}$$

in $L_t^1(\mathbb{R})$. In what follows, we simply write $\widetilde{v}$ as $v$. Since $V$ is real-valued and (A3),

$$V(v(t+h)) - V(v(t)) = 2\mathrm{Re} \left( \overline{(v(t+h) - v(t))} \int_0^1 f(v(t) + \theta(v(t+h) - v(t))) \, d\theta \right)$$

holds. Hence we see from Hölder's inequality that

$$\left\| \int_{\mathbb{R}^n} \frac{V(v(\cdot + h)) - V(v(\cdot))}{h} \, dx - 2 \operatorname{Re} \left\langle \partial_t v(\cdot), \overline{f(v(\cdot))} \right\rangle \right\|_{L^1(\mathbb{R})}$$

$$\lesssim \left| \int_{\mathbb{R}} \left\langle \partial_t v(t), \overline{\int_0^1 f(v(t) + \theta(v(t+h) - v(t))) - f(v(t)) \, d\theta} \right\rangle dt \right|$$

(A.1)
$$+ \left\| \frac{v(\cdot + h) - v(\cdot)}{h} - \partial_t v(\cdot) \right\|_{L_t^\infty H^{-1} \cap L_t^q W^{-1,r}}$$

$$\times \left( \|f(v)\|_{L_t^1 H^1 + L_t^{q'} W^{1,r'}} \right.$$

$$\left. + \int_0^1 \|f(v(t) + \theta(v(\cdot + h) - v(\cdot)))\|_{L_t^1 H^1 + L_t^{q'} W^{1,r'}} d\theta \right),$$

where the first term on the right-hand side of (A.1) makes sense as the duality pairing on $(L_t^1 H^1 + L_t^{q'} W^{1,r'}) \times (L_t^\infty H^{-1} \cap L_t^q W^{-1,r})$. Note that we employ the identity $(f, g) = \langle f, \overline{g} \rangle_{W^{-1,r} \times W^{1,r'}}$ for $f \in L^2 \cap W^{-1,r}$ and $g \in W^{1,r'}$, which follows from the use of $(1 - \varepsilon \Delta)^{-1}$. From $H^1 \hookrightarrow L^r$, it follows that for any $t$,

$$\left\| \int_0^1 f(v(t) + \theta(v(t+h) - v(t))) - f(v(t)) \, d\theta \right\|_{L^2 + L^{r'}}$$
$$\lesssim (1 + \|v(t)\|_{H^1} + \|v(t+h)\|_{H^1})^{p-1} \|v(t+h) - v(t)\|_{H^1} \to 0$$

as $h \to 0$. Similarly to the above, we have

$$\int_0^1 \|f(v(t) + \theta(v(\cdot + h) - v(\cdot)))\|_{L_t^1 H^1 + L_t^{q'} W^{1,r'}} d\theta \lesssim \left( 1 + \|v\|_{L_t^\infty H^1}^{p-1} \right) \|v\|_{L_t^1 H^1 + L_t^{q'} W^{1,r'}}.$$

By the density argument, one can approximate $\partial_t v$ by functions in $C_0^\infty(\mathbb{R} \times \mathbb{R}^n)$, and hence this implies that the first term on the right-hand side of (A.1) converges to zero. To finish the proof, we shall show

$$\left\| \frac{v(\cdot + h) - v(\cdot)}{h} - \partial_t v(\cdot) \right\|_{L_t^\infty H^{-1} \cap L_t^q W^{-1,r}} \to 0$$

as $h \to 0$. Let us only treat the case in $L_t^q W^{-1,r}$, because the case in $L_t^\infty H^{-1}$ easily follows from the strong differentiability of $v$ in $H^{-1}$ and the uniform continuity of $\partial_t v$ in $H^{-1}$. We see from $\partial_t v \in C([-T, T]; H^{-1})$ that

$$v(t+h) - v(t) = \int_t^{t+h} \partial_t v(s) \, ds$$

in $\mathscr{S}'(\mathbb{R}^n)$. Combining the change of variable with Hölder's inequality, one has

$$\left\| \frac{v(\cdot + h) - v(\cdot)}{h} - \partial_t v(\cdot) \right\|_{L^q(\mathbb{R}; W^{-1,r})}$$

$$\leqslant |h|^{\frac{1}{q'} - 1} \left( \int_{\mathbb{R}} \left| \int_t^{t+h} \|\partial_t v(s) - \partial_t v(t)\|_{W^{-1,r}}^q \, ds \right| dt \right)^{\frac{1}{q}}$$

$$= |h|^{\frac{1}{q'} - 1} \left| \int_0^h \left( \int_{\mathbb{R}} \|\partial_t v(t+s) - \partial_t v(t)\|_{W^{-1,r}}^q dt \right) ds \right|^{\frac{1}{q}}$$

$$\leqslant \sup_{|s| \leqslant |h|} \|\partial_t v(\cdot + s) - \partial_t v(\cdot)\|_{L_t^q W^{-1,r}}.$$

Thanks to the continuity of the shift operator in $L_t^q W^{-1,r}$, which is given by the density argument, the above last term converges to zero as $h \to 0$. This completes the proof. $\square$

## Acknowledgments

H.M. was supported by JSPS KAKENHI Grant Number 22K13941. T.O. was supported by JSPS KAKENHI Grant Number 24H00024.

Conflict of Interest: The authors declare that they have no conflict of interest.

Data Availability: Data sharing is not applicable to this article as no datasets were generated or analysed during the current study.

Department of Mathematics, Faculty of Science, SaitamaUniversity, Saitama 338-8570, Japan

Department of Mathematics, Faculty of Science, SaitamaUniversity, Saitama 338-8570, Japan
*Email address*: machihara@rimath.saitama-u.ac.jp

Teacher Training Courses, Faculty of Education, Kagawa University, Takamatsu, gawa 760-8522, Japan
*Email address*: miyazaki.hayato@kagawa-u.ac.jp

Department of Applied Physics, Waseda University, Tokyo 169-8555, Japan
*Email address*: txozawa@waseda.jp